\documentclass[11pt]{amsart}
\usepackage{amscd,amssymb}
\usepackage[all]{xy}
\usepackage[plainpages,backref,urlcolor=blue]{hyperref}

\newtheorem{theorem}{Theorem}[section]
\newtheorem{corollary}[theorem]{Corollary}
\newtheorem{lemma}[theorem]{Lemma}
\newtheorem{prop}[theorem]{Proposition}

\theoremstyle{definition}
\newtheorem{definition}[theorem]{Definition}
\newtheorem{example}[theorem]{Example}
\newtheorem{remark}[theorem]{Remark}
\newtheorem{question}[theorem]{Question}
\newtheorem{problem}[theorem]{Problem}

\newcommand{\Z}{\mathbb{Z}}
\newcommand{\Q}{\mathbb{Q}}
\newcommand{\R}{\mathbb{R}}
\newcommand{\C}{\mathbb{C}}

\newcommand{\T}{\mathbb{T}}
\newcommand{\CP}{\mathbb{CP}}
\newcommand{\RP}{\mathbb{RP}}

\newcommand{\LL}{\mathbb{L}}
\newcommand{\sV}{{\sf V}}
\newcommand{\sE}{{\sf E}}

\renewcommand{\k}{\Bbbk}
\newcommand{\RR}{\mathcal{R}}
\newcommand{\VV}{\mathcal{V}}
\newcommand{\A}{{\mathcal{A}}}

\newcommand{\CC}{{\mathcal{C}}}

\newcommand{\h}{\mathfrak{h}}
\newcommand{\m}{\mathfrak{m}}

\newcommand{\G}{\Gamma}

\DeclareMathOperator{\rank}{rank}
\DeclareMathOperator{\gr}{gr}
\DeclareMathOperator{\im}{im}
\DeclareMathOperator{\coker}{coker}

\DeclareMathOperator{\ideal}{ideal}
\DeclareMathOperator{\id}{id}

\DeclareMathOperator{\Sp}{Sp}

\DeclareMathOperator{\SL}{SL}
\DeclareMathOperator{\SU}{SU}

\DeclareMathOperator{\Hom}{{Hom}}
\DeclareMathOperator{\Tor}{{Tor}}

\DeclareMathOperator{\Lie}{Lie}
\DeclareMathOperator{\tr}{tr}
\DeclareMathOperator{\Aut}{Aut}
\newcommand{\dR}{\scriptscriptstyle{\rm dR}}
\newcommand{\PL}{\scriptscriptstyle{\rm PL}}

\newcommand{\same}{\Longleftrightarrow}
\newcommand{\surj}{\twoheadrightarrow}
\newcommand{\inj}{\hookrightarrow}

\newcommand{\abs}[1]{\left| #1 \right|}

\def\angl#1{{\langle #1\rangle}}
\def\set#1{{\{ #1\}}}

\begin{document}

\title[Geometric and algebraic aspects of $1$-formality]{%
Geometric and algebraic aspects of $1$-formality} 

\author[Stefan Papadima]{Stefan Papadima$^1$}
\address{Institute of Mathematics Simion Stoilow, 
P.O. Box 1-764,
RO-014700 Bucharest, Romania}
\email{Stefan.Papadima@imar.ro}
\thanks{$^1$Partially supported by grant 
CNCSIS ID-1189/2009-2011 of the  
Romanian Ministry of Education and Research}

\author[Alexandru Suciu]{Alexandru~I.~Suciu$^2$}
\address{Department of Mathematics,
Northeastern University,
Boston, MA 02115, USA}
\email{a.suciu@neu.edu}
\urladdr{http://www.math.neu.edu/\~{}suciu}
\thanks{$^2$Partially supported by National Security Agency grant H98230-09-1-0012, 
and an ENHANCE grant from Northeastern University}



\dedicatory{Dedicated to Professor Stere Ianu\c{s} on the 
occasion of his seventieth birthday}

\begin{abstract}
Formality is a topological property, defined in terms  
of Sullivan's model for a space. In the simply-connected 
setting, a space is formal if its rational homotopy type is 
determined by the rational cohomology ring. In the general 
setting, the weaker $1$-formality property allows one to 
reconstruct the rational pro-unipotent completion of the 
fundamental group, solely from the cup products of 
degree $1$ cohomology classes.  

In this note, we survey various facets of formality, with 
emphasis on the geometric and algebraic implications of 
$1$-formality, and its relations to the cohomology jump loci 
and the Bieri--Neumann--Strebel invariant.  We also produce 
examples of $4$-manifolds $W$ such that, for every 
compact K\"{a}hler manifold $M$, the product $M\times W$ 
has the rational homotopy type of a K\"{a}hler manifold, 
yet $M\times W$ admits no K\"{a}hler metric.
\end{abstract}
\maketitle

\tableofcontents

\section{Introduction}
\label{sect:intro}

\subsection{Formality of spaces}
\label{intro:formal}
The question whether one can reconstruct the homotopy type 
of a space from homological data goes back to the beginnings 
of Algebraic Topology.  It was recognized by H.~Poincar\'{e} 
himself that homology is not enough:  for a path-connected 
space $X$, the first homology group, $H_1(X,\Z)$, only records 
the abelianization of the fundamental group, $\pi_1(X)$.  
Even in the simply-connected setting, homology by itself 
fails to detect the Hopf map, $S^3\to S^2$.  On the other 
hand, if one looks at the de~Rham algebra of differential 
forms on $S^n$, one can reconstitute in a purely algebraic 
fashion all the higher homotopy groups of $S^n$, modulo torsion.  

These considerations lead to the notion of {\em formality}\/ 
of a space, as formulated by D.~Sullivan in his foundational work 
on rational homotopy theory \cite{Su77}.  The general definition  
(which we review in \S\ref{sec:dga}) involves a certain commutative 
differential graded algebra, $A_{\PL}(X,\R)$, attached to a 
path-connected space $X$. 
Formality amounts to this cdga being 
related by a chain of quasi-isomorphisms to the cohomology 
algebra, $H^*(X,\R)$, viewed as a cdga with the zero differential.  
In the case when $X$ is a smooth manifold, Sullivan's $A_{\PL}$ 
algebra may be replaced by de~Rham's algebra, leading to the 
following basic principle in rational homotopy theory:
``The manner in which a closed form which is zero in 
cohomology actually becomes exact contains geometric 
information'', cf.~\cite[p.~253]{DGMS}. 

If $X$ is a simply-connected formal space with finite Betti numbers, 
then one can 
build the whole rational Postnikov tower of $X$ (in particular, 
$\pi_*(X)\otimes \Q$), in a purely ``formal" way, just from the 
rational cohomology ring.  On the other hand, one cannot 
hope to do this for an arbitrary formal space of finite type, 
unless $\pi_1(X)$ 
is nilpotent, and acts unipotently on the higher homotopy groups. 
For example, consider the real projective plane, $\RP^2$.  Clearly, 
$\widetilde{H}^*(\RP^2,\Q)=0$; in fact, $\RP^2$ is a formal space, 
with the $\Q$-homotopy type of a point. Yet $\pi_2(\RP^2)=\Z$, 
and so the projective plane does not have the same rational 
homotopy groups as a point. 

\subsection{$1$-Formality of groups}
\label{intro:1formal}
A fruitful way to look at formality in the non-simply-connected 
setting is through the prism of various Lie algebras attached 
to the fundamental group.  If the space $X$ is formal, then  
the group $G=\pi_1(X)$ is {\em $1$-formal}. This means that
the Malcev completion of $G$, as defined by D.~Quillen in \cite{Qu},
is isomorphic, as a filtered Lie algebra, to the completion with 
respect to degree of a quadratic Lie algebra.  In other words, 
the Malcev completion of $G$---and thus, the rational pro-unipotent 
completion of $G$---can be reconstituted from the cup-product 
map, $\cup_{G}\colon H^1(G,\Q)\wedge H^1(G,\Q) \to H^2(G,\Q)$,
more precisely, from the corestriction to its image, $\mu_G$. 
See \S\ref{sec:liealg} for more details.

The main goal of this (non-exhaustive) survey is to present answers
to the following natural question: Given a $1$-formal group $G$, 
what kind of algebraic/topological/geome--
tric information about the 
group can be extracted from the map $\mu_G$?

In the presence of $1$-formality, various other algebraic objects 
attached to $G$ can be recovered from the cup-products in 
degree $1$.  Among these objects, we treat in \S\ref{sec:liealg} 
the graded Lie algebras (modulo torsion) associated to the lower 
central series of the solvable quotients of $G$.  

\subsection{Examples}
\label{intro:examples}
We will come back to the above question in \S\ref{intro:jump}.  
But first, let us address another natural question: 
What sort of conditions insure that a space $X$ is formal---%
or, that a group $G$ is $1$-formal? In essentially all known 
examples, formality follows from one of the following 
reasons:  the cohomology ring of $X$ has some special 
properties; the homotopy groups of $X$ vanish up to a 
certain degree;  $X$ supports some special geometric 
structures; or $G$ admits some distinguished type of presentation. 
Furthermore, there are a number of formality-preserving 
constructions which can be very useful in practice.  

We provide a variety of examples, 
indicating how formality can be deduced from homological, 
geometric, or group-theoretic arguments. We start 
in \S\ref{sec:coho formal} with the connection 
between the algebraic properties of the cohomology ring 
$A=H^*(X,\C)$ and the formality properties of the space $X$. 
In \S\ref{sec:geom}, we discuss the formality properties of 
highly connected manifolds, as well as K\"{a}hler and 
quasi-K\"{a}hler manifolds, smooth affine varieties, and 
Milnor fibers.  Finally, we delineate in \S\ref{sec:gp} several 
classes of $1$-formal groups, including Artin groups, 
Bestvina--Brady groups, and pure welded braid groups. 

\subsection{Cohomology jump loci}
\label{intro:jump}
The $1$-formality property of a group $G$ puts strong restrictions 
on the structure of its cohomology jumping loci.  We survey this 
subject in \S\ref{sec:jump}, where we discuss the characteristic 
varieties $\VV_d(G)$, the resonance varieties $\RR_d(G)$, 
and the relationship between the two.  For a $1$-formal 
group $G$, it turns out that the analytic germs at the 
origin of the characteristic varieties can be reconstructed 
from the  map $\mu_G$.  Furthermore, the analysis of the 
qualitative properties of the resonance varieties of $G$ 
reveals subtle constraints on the cup-product map, imposed 
by $1$-formality.  As we explain in Example \ref{ex=euclid}, 
this is a striking phenomenon, peculiar to 
{\em non-simply-connected}\/ rational homotopy theory.

We also discuss in \S\ref{sec:jump} the Bieri--Neumann--Strebel 
invariant, $\Sigma^1(G)$---a rather enigmatic object, which 
controls the finiteness properties of normal subgroups of $G$ with
abelian quotient. Again, the $1$-formality assumption plays 
a significant role, and yields an upper bound for the BNS 
invariant $\Sigma^1(G)$, depending solely on $\mu_G$.

Under suitable geometric assumptions, the restrictions imposed 
by $1$-formality on the cohomology jumping loci are strong enough 
to lead to complete classification results.   In \S\ref{sect:serre}, 
we describe how cup-products in degree $1$ can detect the 
realizability of a $1$-formal group $G$, as the fundamental 
group of a (quasi-) K\"{a}hler manifold.  This method applies 
to a wide range of groups, including right-angled Artin groups, 
Bestvina--Brady groups, and $3$-manifold groups. 

\subsection{Algebraic monodromy}
\label{intro:mono}
Given a locally trivial, smooth fibration, $F\inj M\to B$, there 
is an associated algebraic monodromy action of 
$\pi_1(B)$ on $H_*(F)$. Likewise, every group extension, 
$N\inj G\surj Q$, gives rise to a monodromy action of 
$Q$ on $H_*(N)$, induced by conjugation in $G$. 
In \S\ref{sect:mono}, we examine the interplay between 
monodromy and $1$-formality.  For Artin kernels, 
triviality of the monodromy action insures $1$-formality, 
whereas for fibrations over the circle, formality of the 
total space implies absence of monodromy Jordan blocks of size 
greater than $1$ for the eigenvalue $1$. 

Needless to say, the formal theory has its limitations. We 
illustrate this point in the last section (where all the new 
material is concentrated), with a family of examples 
originating from symplectic geometry.  The upshot is 
the following theorem, which we prove in 
\S\ref{sect:epilogue}, by analyzing the algebraic 
monodromy of $3$-dimensional mapping tori.

\begin{theorem}
\label{thm:new}
There exist infinitely many closed, orientable, formal 
$4$-manifolds $W$ such that, for every compact 
K\"{a}hler manifold $M$, the following hold. 
\begin{enumerate}
\item The manifold $M\times W$ has the same rational homotopy 
type as the K\"{a}hler manifold $M\times T^2 \times S^2$. 
\item The manifold $M\times W$ admits no K\"{a}hler metric.  
\end{enumerate}
\end{theorem}

\subsection{Conventions}
\label{intro:conv}
Except otherwise stated, by a space we always mean a 
topological space having the homotopy type of a connected 
polyhedron with finite $1$-skeleton.  Similarly, all groups 
we consider here are assumed to be finitely generated. 
The typical examples we have in mind are compact, 
connected, smooth manifolds and their fundamental groups. 
We say that a manifold is {\em closed}\/ if it is smooth, 
compact, connected, and boundaryless.  Coefficients are 
usually taken in a field $\k$ of characteristic zero; when 
coefficients are not mentioned, the default is $\k=\C$. 

\section{From spaces to differential graded algebras}
\label{sec:dga}

We start with Sullivan's construction of an algebraic model  
encoding the rational homotopy type of a space. 

\subsection{Differential graded algebras}
\label{subsec:dga}
Fix a ground field $\k$, of characteristic $0$.  
A {\em commutative differential graded algebra} (for short, 
a cdga), is a graded $\k$-algebra $A$, endowed with 
a differential $d_A\colon A\to A$ of degree $1$. We 
assume here commutativity in the graded sense, that is, 
$ab =(-1)^{\abs{a}\abs{b}} ba$, for every homogeneous 
elements $a, b\in A$, where $\abs{a}$ denotes the degree of $a$. 

A cdga morphism is a {\em quasi-isomorphism}\/ if it induces 
an isomorphism in cohomology.  Two commutative differential 
graded algebras, $A$ and $B$, are said to be {\em weakly 
equivalent}\/ if there is a zig-zag of quasi-isomorphisms 
(going both ways), connecting $A$ to $B$. 

\begin{definition}
\label{def:formal}
A cdga is {\em formal}\/ if it is weakly equivalent to its cohomology 
algebra, endowed with the zero differential. 
\end{definition}

We will also consider the following natural notion of ``partial" formality, up 
to some degree $q\ge 1$. 

\begin{definition}
\label{def:qformal}
A cdga $(A,d_A)$ is {\em $q$-formal}\/ if there is a zig-zag 
of morphisms connecting $(A,d_A)$ to $(H^*(A, d_A), d=0)$, 
with each one of these maps inducing an isomorphism in cohomology 
up to degree $q$, and a monomorphism in degree $q+1$.  
\end{definition}

\subsection{Models of spaces}
\label{subsec:models}

Let $X$ be a space.  (Recall we are tacitly assuming $X$ 
is homotopy equivalent to a connected polyhedron with 
finite $1$-skeleton.)  In \cite{Su77}, Sullivan constructs 
an algebra $A_{\PL}(X, \k)$ of polynomial differential forms 
on $X$ with coefficients in $\k$, and provides it with a 
natural cdga structure.  A {\em model}\/ for $X$, over 
the field $\k$, is a cdga weakly equivalent 
to $A^*_{\PL}(X, \k)$.  Two spaces, $X$ and $Y$, have 
the same {\em $\k$-homotopy type}\/ if $A^*_{\PL}(X, \k)$ 
and $A^*_{\PL}(Y, \k)$ are weakly equivalent. 

The space $X$ is said to be formal (over $\k$) if Sullivan's 
algebra $A_{\PL}(X, \k)$ is formal.  Likewise, $X$ is $q$-formal 
if this cdga is $q$-formal.  When $X$ is a smooth manifold, 
and $\k=\R$ or $\C$, we may replace in the definition the 
algebra of polynomial forms by the corresponding de~Rham 
algebra of differential forms, $\Omega_{\dR}(X,\k)$; 
see for instance \cite{FHT}. 

Clearly, formality implies partial formality, and $q$-formality implies 
$r$-formality, for all $r\le q$. We refer to M\u{a}cinic \cite{Mac} for a
study of $q$-formality in the range $q\ge 2$.  Here, we are 
primarily interested in $1$-formality. 

\begin{remark}
\label{rem:kg1}
The $1$-formality property of a space $X$ depends only on 
its fundamental group, $G=\pi_1(X)$. Indeed, let 
$f\colon X\to K(G,1)$ be a classifying map.  Then, 
the induced homomorphism, $f^*\colon H^i(G,\k) \to H^i(X,\k)$, 
is an isomorphism for $i=1$ and a monomorphism for $i=2$. 
The claim follows. 
\end{remark}

\subsection{Massey products}
\label{ss=massey}

The first approach to $1$-formality is in terms of certain higher order 
structures from the algebraic homotopy theory of differential graded 
algebras, called {\em Massey products}. This point of view has been 
extensively used to compare symplectic and K\"{a}hler structures on 
manifolds, see for instance \cite{FOT, FGM, TO}. 

Since we will not pursue this direction here, we will avoid 
precise definitions, and simply say that ``a group $G$ 
is $1$-formal if and only if all Massey products of 
elements from $H^1(G,\k)$ vanish uniformly, for length 
at least $3$''; compare with \cite[p.~262]{DGMS}.  
(The length $2$ Massey products are simply the cup-products.)

For a differential graded algebra $(A,d)$, the Massey 
product of classes $\alpha_1,\alpha_2,\alpha_3\in H^1(A)$ 
is defined, provided $\alpha_1 \alpha_2 =\alpha_2\alpha_3=0$.  
Pick representative cocycles $a_i$ for $\alpha_i$, 
and elements $y,z\in A$ such that $dy=a_1a_2$ 
and $dz=a_2a_3$. It is readily seen that $ya_3+a_1z$ 
is a cocyle.  The set of cohomology classes of all 
such cocycles is the Massey triple product
$\angl{\alpha_1,\alpha_2,\alpha_3}$.   
The image of this set  in the quotient ring 
$H^*(A)/( \alpha_1,\alpha_3)$ is a well-defined 
element in degree $2$; we say 
$\angl{\alpha_1,\alpha_2,\alpha_3}$ 
is {\em non-vanishing}\/ if this element is not $0$.  

\begin{example}
\label{ex:triple}
Let $M=G_{\R}/G_{\Z}$ be the $3$-dimensional Heisenberg 
nilmanifold, where $G_{\R}$ is the group of real, 
unipotent $3\times 3$ matrices, and $G_{\Z}=\pi_1(M)$ 
is the subgroup of integral matrices in $G_{\R}$.  One may 
use invariant forms to obtain the following simple model 
$(A,d)$ for the Heisenberg manifold: $A=\bigwedge (a,b,z)$ 
is the exterior algebra on the indicated generators in degree 
$1$, and the differential is given by $da=db=0$, $dz=ab$. 
Clearly, $\cup_{M} =0$, and $\angl{[a],[a],[b]}= [az]$, 
with trivial indeterminacy.  Since $[az]\ne 0$, 
the manifold $M$ is not $1$-formal. 
\end{example}

\section{From groups to Lie algebras}
\label{sec:liealg}

There is a dual approach to $1$-formality, based on Lie algebras. 
In this section, we review the construction of several Lie algebras 
associated to a group, and how these Lie algebras are related 
to each other, in the case when the group is $1$-formal.  

\subsection{Holonomy Lie algebras}
\label{subsec:holo}
Given a space $X$, let 
$\cup_X \colon H^1(X,\k)\wedge H^1(X,\k)\to H^2(X,\k)$
be the cup-product map in degree one, and let 
$\partial_X\colon H_2(X, \k)\to  H_1(X, \k) \wedge H_1(X, \k)$ 
be the comultiplication map.  

\begin{definition}[K.-T. Chen \cite{Ch}]
\label{def:holo lie}
The {\em holonomy Lie algebra}\/ of $X$ over $\k$ is the 
quotient
\begin{equation*}
\label{eq:hlie}
\h (X,\k) = \Lie(H_1(X, \k))/ \ideal (\im (\partial_X)),
\end{equation*}
where $\Lie(H_1(X, \k))$ is the free Lie algebra over $\k$, 
generated  by $H_1(X, \k)$ and graded by bracket length.  
\end{definition}

Note that the defining ideal of $\h(X,\k)$ is a homogeneous (in fact, 
quadratic) ideal.  Hence, the holonomy Lie algebra inherits 
a natural grading from the free Lie algebra.  Further, note 
that $\h (X,\k)$ depends only on $\mu_X$, the corestriction 
of $\cup_X$ to its image.  

Now assume $X$ is a path-connected space, and let 
$G=\pi_1(X)$. Define the holonomy Lie algebra 
of the group $G$ as that of a classifying space $K(G,1)$:
\[
\h(G,\k):=\h (K(G,1),\k).
\]
It is readily seen that $\mu_X =\mu_G $.  It follows that 
$\h(X,\k)=\h(G,\k)$.

\subsection{Malcev Lie algebras}
\label{subsec:malcev}
Next, we recall some notions from  \cite[Appendix~A]{Qu}.  
A {\em Malcev Lie algebra}\/ is a Lie algebra over $\k$, 
endowed with a decreasing, complete $\k$-vector space 
filtration, satisfying certain axioms (see \S\ref{subsec:lcs} 
below for more details). For example, the completion of 
$\h(X,\k)$ with respect to the degree filtration is a Malcev 
Lie algebra, denoted $\widehat{\h}(X,\k)$.  

In \cite{Qu}, Quillen associates to every group $G$, in a functorial 
way, a Malcev Lie algebra, denoted $\m(G,\k)$.  This object, 
called the {\em Malcev completion}\/ of $G$, captures 
the properties of the torsion-free nilpotent quotients of $G$. 

Here is a concrete way to describe it. 
The group algebra $\k{G}$ has a natural Hopf algebra structure,
with comultiplication given by $\Delta(g)=g\otimes g$, and
counit the augmentation map. Let $I$ be the augmentation ideal. 
One verifies that the Hopf algebra structure on $\k{G}$ extends 
to the $I$-adic completion, $\widehat{\k{G}}=\varprojlim_r \k{G}/I^r$.  
Finally,  $\m(G,\k)$ coincides with the Lie algebra of primitive 
elements in $\widehat{\k{G}}$, endowed with the inverse limit 
filtration.

\begin{theorem}[Sullivan \cite{Su77}]
\label{thm=1fmal}
A group $G$ is $1$-formal, over $\k$, if and only if 
$\m(G,\k)\cong \widehat{\h}(G,\k)$, as filtered Lie algebras.
\end{theorem}

Consequently, if $G$ is a $1$-formal group, then the co-restriction to 
the image of the cup-product map, $\mu_G$,  determines the 
Malcev completion $\m(G)$.  

The next result is folklore.  A proof is given in \cite[Lemma 2.9]{DPS-serre} 
for finitely presented groups, but the argument given there 
works as well for finitely generated groups.

\begin{lemma}
\label{lem:quadratic}
A group $G$ is $1$-formal, over $\k$, if and only if the 
Malcev Lie algebra $\m(G,\k)$ is isomorphic, as a filtered Lie algebra, 
to the completion with respect to degree of a quadratic Lie algebra.
\end{lemma}

\begin{example}
\label{ex:free}
Let $F_n$ be the free group of rank $n\ge 0$.  Clearly, $H_1(F_n,\k)=\k^n$ 
and $H_2(F_n,\k)=0$; in particular, $\cup_{F_n}=0$.  Thus, $\h(F_n,\k)$ 
is isomorphic to $\LL_n=\Lie(\k^n)$, the free Lie algebra of rank $n$, 
over $\k$.  It is readily checked that $\m(F_n,\k)=\widehat{\LL}_n$.  
It follows from Theorem \ref{thm=1fmal}  (or Lemma \ref{lem:quadratic}) 
that $F_n$  is $1$-formal. 
\end{example}

\begin{example}
\label{ex:surface}
Let $\Sigma_g$ be the Riemann surface of genus $g\ge 1$.   
The group $G=\pi_1(\Sigma_g)$ is generated 
by $x_1,y_1,\dots,x_g,y_g$, subject to the single relation 
$[x_1,y_1]\cdots [x_g,y_g]=1$, where 
$[x,y]=xyx^{-1}y^{-1}$ is the group commutator. It is readily checked that $\h(G,\k)$ 
is the quotient of the free Lie algebra on $x_1,y_1,\dots,x_g,y_g$  
by the ideal generated by $[x_1,y_1]+ \cdots + [x_g,y_g]$. 
A further computation shows that $\m(G,\k)\cong \widehat{\h}(G,\k)$, 
and thus, $G$ is $1$-formal. 
\end{example}

\subsection{Associated graded Lie algebra}
\label{subsec:lcs}

Given an arbitrary group $G$, the lower central series (for 
short, LCS) of $G$  is defined inductively by $\gamma_1 G=G$ 
and $\gamma_{k+1}G =[\gamma_k G,G]$. 
The {\em associated graded Lie algebra}, $\gr(G)$, 
is the direct sum of the successive LCS quotients, 
\begin{equation*}
\label{eq:grg}
\gr(G)= \bigoplus\nolimits_{k\ge 1} \gamma_k G/ \gamma_{k+1} G,
\end{equation*} 
with Lie bracket induced from the group commutator. 

The Malcev filtration $\{ F_s\}_{s\ge 1}$ on $\m(G,\k)$ 
is required to satisfy $[F_s,F_t]\subseteq F_{s+t}$, for all $s,t$. 
Consequently, the associated graded vector space, 
$\gr (\m (G,\k))= \bigoplus_{s\ge 1} F_s/F_{s+1}$,
inherits a natural Lie algebra structure, compatible 
with the grading. The basic property of the Malcev 
completion is that 
\begin{equation*}
\label{eq:grmalcev}
\gr (G)\otimes \k\cong \gr (\m (G,\k)),
\end{equation*} 
as Lie algebras with grading. We infer from 
Theorem \ref{thm=1fmal} that $\mu_G$ also 
determines $\gr(G)$ modulo torsion, in the $1$-formal case.

\begin{corollary}
\label{cor=1fgr}
If the group $G$ is $1$-formal, then $\gr (G)\otimes \k\cong \h(G,\k)$, 
as Lie algebras with grading.
\end{corollary}

The {\em derived series} of a group $G$ is defined inductively by 
$G^{(0)}=G$ and $G^{(i+1)}=[G^{(i)},G^{(i)}]$. The derived series 
of a Lie algebra is constructed similarly. The next result extends 
Corollary \ref{cor=1fgr} to solvable quotients. 

\begin{theorem}[\cite{PS-chenlie}]
\label{thm=solvq}
If the group $G$ is $1$-formal, then for each $i\ge 1$, 
\[
\gr (G/G^{(i)})\otimes \k\cong \h(G,\k)/\h^{(i)}(G,\k), 
\]
as Lie algebras with grading.
\end{theorem}

The proof, given in \cite[Theorem 4.2]{PS-chenlie} in the case 
when $G$ is finitely presented, works as well when $G$ 
is finitely generated.

\section{Cohomology ring and formality}
\label{sec:coho formal}

We now turn to the task of delineating a set of conditions, 
sufficient to guarantee the formality of a space, or the $1$-formality 
of a group. We start in this section with some cohomological considerations.

\subsection{Small first Betti number} 
\label{subsec:betti}

As usual, let $G$ be a finitely generated group, and let $\k$ be a field 
of characteristic $0$. Denote by $b_i(G)=\dim_{\k} H_i(G,\k)$ 
the $i$-th Betti number of $G$.  Let us start with a well-known 
(and easy to prove) fact.  

\begin{prop}
\label{prop:b1}
If $b_1(G)\le 1$, then $G$ is $1$-formal. 
\end{prop}

\begin{proof}
When $b_1(G)=0$, the homomorphism $f\colon G\to \set{1}$  
induces an isomorphism on $H^1$ and a monomorphism on $H^2$. 
Likewise, when $b_1(G)=1$, we may pick a homomorphism 
$f\colon G\to \Z$ inducing an isomorphism on $H^1$; it follows 
then that $f$ is injective on $H^2$. Our claim follows from the 
$1$-formality of free groups, discussed in Example \ref{ex:free}.
\end{proof}

This result is optimal.  Indeed, the Heisenberg group $G=G_{\Z}$ 
from Example \ref{ex:triple} has $b_1(G)=2$, and is not $1$-formal.  

In the case of $3$-manifolds, the above proposition was 
sharpened by Fern\'{a}ndez and Mu\~{n}oz, as follows. 

\begin{theorem}[\cite{FM05}]
\label{thm:fm}
Let $M$ be a closed, orientable $3$-manifold with 
$b_1(M)\le 1$.  Then $M$ is formal, and has the same
$\k$-homotopy type as $S^3$ or $S^1\times S^2$.
\end{theorem}

\subsection{Regular sequences} 
\label{subsec:reg sec}

Let $A$ be a graded, graded-commutative algebra (for short, cga) 
over $\k$. A sequence $r_1,\dots, r_n$ of elements of $A$ is 
said to be a {\em regular sequence} if $r_i$ is not a zero-divisor in
$A/(r_1,\dots, r_{i-1})$, for each $i\le n$. 

\begin{theorem}[Sullivan \cite{Su77}]
\label{thm=reg}
If $H^*(M,\k)$ is the quotient of a free cga by an ideal generated 
by a regular sequence, then $M$ is a formal space. In particular, 
freeness of $H^*(M,\k)$ implies formality of $M$.
\end{theorem}

This result provides a large supply of formal spaces, such as
\begin{itemize}
\item rational cohomology spheres; 
\item rational cohomology tori; 
\item compact connected Lie groups $G$, as well as their 
classifying spaces, $BG$; 
\item homogeneous spaces of the form $G/K$, with 
$\rank G=\rank K$;  
\item Eilenberg-MacLane spaces $K(\pi, n)$ with $n\ge 2$.
\end{itemize}

In particular, if $X$ is the complement of a knotted sphere 
in $S^{n}$, $n\ge 3$, then $X$ is a formal space. 

On the other hand, not all homogeneous spaces are formal:  
for instance, $\Sp(5)/\SU(5)$ is not, see e.g.~\cite[p.~\!143]{FOT}.  
And Eilenberg-MacLane spaces $K(\pi, 1)$ need not be formal: 
for example, if $\pi$ is a torsion-free, finitely generated 
nilpotent group, then $K(\pi, 1)$ is formal if and only if 
$\pi$ is abelian, see e.g.~\cite[p.~\!120]{FOT}. 

\subsection{From partial to full formality} 
\label{subsec:koszul}

In general, $q$-formality is strictly weaker than formality, for $q<\infty$. 
Nevertheless, under favorable circumstances the two notions are equivalent. 
One result of this type is related to the well-known Koszul property 
from homological algebra. 

Let $A$ be a graded $\k$-algebra with $A^0=\k$.  Then, 
\[
\Tor^A (\k ,\k)= \bigoplus_{s,t\ge 0} \Tor^A_s (\k ,\k)_t
\]
is a bigraded vector space: the index $s$ denotes the usual 
homological degree, while $t$ stands for the internal degree, 
coming from the grading of $A^*$. The algebra $A$ is said 
to be a {\em Koszul algebra}\/ if $\Tor^A_s (\k ,\k)_t=0$, for $s\ne t$.   

\begin{theorem}[\cite{PY, PS-artin}]
\label{thm=koszul}
Let $X$ be a connected CW-complex with finite skeleta. Suppose 
$H^*(X,\k)$ is a Koszul algebra.  Then $X$ is $1$-formal if and 
only if $X$ is formal. 
\end{theorem}

The following result of M\u{a}cinic ties partial formality to full 
formality, under a completely different homological assumption.

\begin{theorem}[\cite{Mac}]
\label{thm=qfvsf}
Let $X$ be a space with the property that $H^i(X,\k)=0$, for $i>q+1$. 
If $X$ is $q$-formal, then $X$ is formal.
\end{theorem}

Here is an immediate consequence. 

\begin{corollary}
\label{cor:qff}
Every $q$-formal CW-complex of dimension at most $q+1$ is 
formal.
\end{corollary}

In particular, if $G$ is a finitely-generated, $1$-formal group, 
and $K$ is a $2$-complex with $\pi_1(K)=G$, then $K$ 
is formal.  See \cite[Lemma 2.10]{FM07} for a
result similar to Corollary \ref{cor:qff}, involving 
a different (more restrictive) notion of $q$-formality,
introduced by Fern\'andez and Mu\~noz in \cite{FM07}.

\section{Manifolds and geometric structures}
\label{sec:geom}

In this section, we look at some of the ways in which formality 
of a space is implied by relevant topological properties 
or geometric structures.

\subsection{Cell complexes and manifolds} 
\label{subsec:manifolds}
We start with connectivity properties: roughly speaking, the 
more highly connected a finite-dimensional CW-complex is, 
the more likely it is to be 
formal.  This was made precise by Stasheff \cite{Sta83}, 
as follows.  Let $X$ be a $k$-connected CW-complex of 
dimension $n$; if $n \le 3k+1$, then $X$ is formal.  This 
is the best possible bound:  attaching a cell $e^{3k+2}$ to the 
wedge $S^{k+1}\vee S^{k+1}$ via the iterated Whitehead 
product $[\iota_1,[\iota_1,\iota_2]]$ yields a non-formal 
CW-complex. 

Formality is preserved under several standard operations 
on (based) CW-complexes with finite Betti numbers.  
For instance, if $X$ and $Y$ 
are formal, then so is the product $X\times Y$ and the wedge 
$X\vee Y$; moreover, a retract of a formal space is formal. 
We refer to \cite{FHT, FOT} for details. 

For closed manifolds, the above dimension bound 
can be relaxed, by using Poincar\'{e} duality.  
As shown by Miller \cite{Mi79}, if $M$ is a closed, 
$k$-connected manifold ($k\ge 0$) of dimension $n \le 4k+2$, 
then $M$ is formal.  In particular, all simply-connected closed
manifolds of dimension at most $6$ are formal.  
Again, this is best possible: as shown by 
Fern\'andez and Mu\~noz in \cite{FM04}, there exist  
closed, simply-connected, non-formal manifolds of 
dimension $7$.  On the other hand, if $M$ is a closed, 
orientable, $k$-connected $n$-manifold with $b_{k+1}(M)=1$, 
then the bound insuring formality can be improved to $n \le 4k+4$, 
see Cavalcanti \cite{Ca}.  

Formality behaves well with respect to certain operations 
on manifolds.  For instance, Stasheff \cite{Sta83} proved the 
following: If $M$ is a closed, simply-connected manifold such 
that $M\setminus \set{*}$ is formal, then $M$ is formal.  
Moreover, if $M$ and $N$ are closed, orientable,
formal manifolds, so is their connected sum, $M\# N$;
see  \cite{FHT}.

\subsection{K\"{a}hler manifolds and smooth algebraic varieties}
\label{subsec:kahler}
Certain geometric structures influence favorably the formality 
of manifolds.   For instance, on a compact Riemannian symmetric 
space $M$, the product of harmonic forms is again harmonic. 
Hodge theory, then, implies the formality of $M$ \cite{Su77}. 

Hodge theory also has strong implications on the topology 
of compact K\"{a}hler manifolds.  If $M$ is such a manifold, 
let $d$ be the exterior derivative, $J$ the complex structure, 
and $d^c=J^{-1}dJ$.  In \cite{DGMS}, Deligne, Griffiths, 
Morgan, and Sullivan showed that the following ``$dd^c$ Lemma" 
holds for $M$: if $\alpha$ is a form which is closed for both $d$ and 
$d^c$, and exact for either $d$ or $d^c$, then $\alpha$ is exact  
for $dd^c$.  Formality ensues:

\begin{theorem}[\cite{DGMS}]
\label{thm:dgms}
All compact K\"{a}hler manifolds are formal.
\end{theorem}

A manifold $M$ is said to be a {\em quasi-K\"{a}hler manifold} if 
$M=\overline{M}\setminus D$, where $\overline{M}$ is 
a compact K\"{a}hler manifold, and $D$ is a normal 
crossing divisor.  For example, smooth, irreducible, 
quasi-projective complex varieties are quasi-K\"{a}hler. 
In \cite[Corollary 10.3]{M}, Morgan establishes the following result. 

\begin{theorem}[\cite{M}]
\label{thm:morgan}
Let $M$ be a smooth, quasi-projective variety.  If the Deligne 
weight filtration space $W_1H^1(M,\C)$ vanishes, then 
$M$ is $1$-formal. 
\end{theorem}

This happens, for instance, when $M$ is the complement of a 
hypersurface in $\CP^n$.  

\begin{example}
\label{ex:curves}
Let $\mathcal{C}$ be an algebraic curve in $\CP^2$.  By the above,  
the complement $M=\CP^2\setminus \mathcal{C}$ is $1$-formal. 
On the other hand, $M$ has the homotopy type of a finite 
CW-complex of dimension $2$.   Thus, by Corollary \ref{cor:qff} 
to M\u{a}cinic's theorem, $M$ is formal.  For a different proof 
of this result, see \cite[Theorem 6.4]{CM}. 
\end{example}

Also note that every smooth, irreducible complex curve $C$ is formal.  
In the compact case, formality follows from the K\"{a}hler property, 
while in the non-compact case, formality follows from the fact that 
$C$ is homotopy equivalent to a (finite) wedge of circles. 

\subsection{Affine varieties and Milnor fibrations}
\label{subsec:affine}
In contrast, smooth, irreducible {\em affine}\/ varieties need 
not be $1$-formal.  A general construction illustrating this phenomenon  
is given in \cite[Proposition 7.2]{DPS-qk}.  Here is a concrete example, 
taken from  \cite{DPS-qk}. 

\begin{example}
\label{ex:affine}
Consider the polynomials $g=x^3+y^3+z^3$ and $f=x+y^2+z^3$.  
Then $M=V(g) \setminus V(f)$ is a smooth affine subvariety of $\C^4$, 
yet $M$ is not $1$-formal. 
\end{example}

An important construction in singularity theory is that of the 
{\em Milnor fibration}.  In its simplest incarnation, this goes 
as follows.  Let $f\in \C[z_0,\dots, z_n]$ be a homogeneous 
polynomial. The restriction $f\colon \C^{n+1} \setminus V(f) \to \C^*$ 
turns out to be a smooth bundle projection.  Clearly, the fiber 
of this bundle, $F(f):=f^{-1}(1)$, is a smooth affine variety, 
having the homotopy type of an $n$-dimensional finite 
CW-complex. When the singularity $(V(f), 0)$ is reduced, 
the Milnor fiber $F(f)$ is connected. The above considerations 
naturally lead to the following question. 

\begin{question}
\label{q=milnorf}
Is the Milnor fiber of a reduced polynomial, $F(f)$, 
always  $1$-formal?
\end{question}

When $f=f_1\cdots f_d$ completely splits as a product of distinct linear 
factors, the reduced variety $V(f)$ is a finite union of hyperplanes in 
$\C^{n+1}$.  In this case, the complement  $M=\C^{n+1} \setminus V(f)$ 
is a formal space. This follows from work by Brieskorn \cite{Br}, who 
showed that the inclusion of the subalgebra generated by the closed 
logarithmic $1$-forms ${df_i}/{f_i}$ into $\Omega^*_{\dR}(M)$
induces an isomorphism in cohomology. We mention that
the above question is open even for hyperplane arrangements.

\section{Group presentations and $1$-formality}
\label{sec:gp}

In this section, we discuss the $1$-formality property 
of several classes of groups, as well as the behavior 
of this property under standard group-theoretic constructions. 

\subsection{Commutator relators and vanishing cup products}
\label{subsec:comm rels}
Suppose $G$ is a finitely presented group.  Then, 
as shown in \cite{P1}, one may find a finite presentation 
for the corresponding Malcev Lie algebra, $\m(G)$. 
In favorable cases, this latter presentation is quadratic, 
and thus one may apply Lemma \ref{lem:quadratic} 
to conclude that $G$ is $1$-formal. 

As we saw in Example \ref{ex:free}, the free group $F_n$ 
has vanishing cup-product map $\cup_{F_n}$, and is 
$1$-formal.   Here is a partial converse.

\begin{prop}[\cite{DPS-jump}]
\label{prop:commrel}
Let $G$ be a group admitting a finite presentation with 
only commutator relators.  If $G$ is $1$-formal and 
$\cup_G=0$, then $G$ is a free group. 
\end{prop}

The proposition shows that, in some sense ``generically," 
$1$-formality does not hold, at least not among 
commutator-relators groups. 

\begin{remark}
\label{rem:heispres}
The above result explains---without resorting to Massey products---why 
the Heisenberg group $G=G_{\Z}$ from Example \ref{ex:triple} cannot 
be $1$-formal.  Indeed, the group has commutator-relators 
presentation $G=\angl{ x, y \mid [x,[x,y]], [y,[x,y]]}$, and $\cup_G=0$. 
\end{remark}

\begin{remark}
\label{rem:purelinks}
Let $L$ be a link in $S^3$ obtained by closing up 
a pure braid, and let $M=S^3\setminus L$ be its complement. 
It is readily seen that the group $G=\pi_1(M)$ admits a 
commutator-relations presentation. Assuming $\cup_G =0$, 
it follows from Proposition \ref{prop:commrel} that $G$ is 
$1$-formal if and only if $L$ is a trivial link.  Therefore, 
we cannot expect to obtain much information from the 
vanishing of $\cup_G$, except, of course, the vanishing 
of all linking numbers. 
\end{remark}

\subsection{Products, coproducts, and extensions}
\label{subsec:operations}
The $1$-formality property behaves well with respect to (finite) 
direct products and coproducts.  

\begin{theorem}[\cite{DPS-jump}]
\label{thm:prod coprod}
Let $G_1$ and $G_2$ be finitely presented, $1$-formal groups. 
Then $G_1\times G_2$ and $G_1*G_2$ are also 
$1$-formal. 
\end{theorem}

By contrast, $1$-formality behaves rather badly with respect to quotients 
and subgroups. Of course, any group is a quotient of a free group,
and thus free groups of rank at least $2$ possess plenty of 
non-$1$-formal quotients. The next example shows that (finitely 
generated) subgroups of $1$-formal groups need not be $1$-formal. 

\begin{example}
\label{ex=subgr}
Let $G=G_{\Z}$ be the Heisenberg group, with presentation 
as in Remark \ref{rem:heispres}. Consider the semidirect 
product $H=G\rtimes_{\varphi} \Z$, defined by the 
automorphism $\varphi\colon G\to G$ given by 
$x\mapsto y$, $y\mapsto xy$. Since clearly 
$b_1(H)=1$, the group $H$ is $1$-formal. Yet, the 
normal subgroup $G=H^{(1)}$ is not $1$-formal.
\end{example}

\subsection{Artin groups}
\label{subsec:artin}

Let $\Gamma=(\sV,\sE,\ell)$ be a labeled finite simplicial graph, 
with vertex set $\sV$, edge set $\sE$, and labeling function 
$\ell \colon \sE\to \Z _{\geq 2}$.  The corresponding  
{\em Artin group} has one generator for each vertex 
$v \in \sV$ and one defining relation 
\[
\underbrace{vwv\cdots}_{\ell(e)} =\underbrace{wvw\cdots}_{\ell(e)}
\]
for each edge $e=\{v,w\}$ in $\sE$.  For example, if 
$\Gamma=K_{n-1}$ is the complete graph on vertices 
$1$ through $n-1$, with labels $\ell(\{i,j\})=2$ if $\abs{i-j}>1$ and  
$\ell(\{i,j\})=3$ if $\abs{i-j}=1$, then the corresponding Artin group 
is the braid group on $n$ strings, $B_n$. 

If $\Gamma=(\sV,\sE)$ is unlabeled, then $G_{\Gamma}$ is 
called a {\it right-angled Artin group}, and is defined by 
commutation relations $vw=wv$, one for each edge 
$\{v,w\} \in\sE$.  These groups interpolate between $\Z^n$ 
(for $\Gamma=K_n$) and $F_n$ (for $\Gamma=\overline{K}_n$), 
and behave nicely with respect to the join operation for graphs: 
$G_{\Gamma*\Gamma'}=G_{\Gamma}\times G_{\Gamma'}$. 

Using the defining presentations, Kapovich and Millson proved 
the following theorem.

\begin{theorem}[\cite{KM}]
\label{thm=km}
All  Artin groups are $1$-formal.
\end{theorem}

By combining Theorems \ref{thm=km} and \ref{thm=koszul}, 
we showed in \cite{PS-artin} that the classifying spaces of 
right-angled Artin groups are formal spaces.  More generally, 
Notbohm and Ray showed that, for each finite simplicial 
complex $L$, the corresponding toric complex $T_L$ is 
formal; see \cite[Remark 5.7]{NR}. 

\subsection{Bestvina--Brady groups}
\label{subsec:bb}
Given a finite simple graph $\G$, let 
\[
N_{\G}=\ker(\nu\colon G_{\G}\surj \Z)
\]
be the kernel of the ``diagonal" epimorphism, sending 
each generator $v$ to $1$.  As shown by Bestvina and 
Brady in \cite{BB}, the homological finiteness properties of the 
group $N_\G$ are intimately connected to the topology of 
the flag complex $\Delta_\G$, that is, the maximal simplicial 
complex with $1$-skeleton $\G$.  For example, $N_{\G}$ 
is finitely generated if and only if $\G$ is connected; and 
$N_{\G}$ is finitely presented if and only if $\Delta_\G$ 
is simply-connected.  

Using the presentation of $N_{\G}$ derived by Dicks and Leary 
in \cite{DL}, we proved in \cite{PS-bb} the following 
result. 

\begin{theorem}[\cite{PS-bb}]
\label{thm:formalbb}
All finitely presented Bestvina--Brady groups $N_{\Gamma}$ 
are $1$-formal.
\end{theorem}

\subsection{Welded braid groups}
\label{subsec:weld}
Let $F_n=\angl{x_1,\dots,x_n}$ be the free group of rank $n\ge 1$. 
The {\em welded braid group}\/ on $n$ strands is the subgroup 
of $\Aut (F_n)$ consisting of those group automorphisms which 
send each generator $x_i$ to a conjugate of another generator. 
The elements for which the associated permutation of the generators 
is the identity form the {\em pure}\/ welded braid subgroup. 

Making use of the presentation given by McCool in \cite{MC}, the 
following result was obtained by Berceanu and Papadima.

\begin{theorem}[\cite{BP}] 
\label{thm=bp}
Pure welded braid groups are $1$-formal.
\end{theorem}

This theorem extends the well-known $1$-formality of pure 
braid groups (which are fundamental groups of hyperplane 
arrangement complements), and leads to an explicit construction 
of the Kontsevich integral for the welded braid groups, 
described in \cite{BP}.

\section{Cohomology jump loci and the BNS invariant}
\label{sec:jump}

In this section, we describe some implications of 
the $1$-formality property on the structure of the 
cohomology jump loci and the Bieri--Neumann--Strebel 
invariant of a space (or a group).

\subsection{Jump loci}
\label{subsec:jumps}

As before, let $X$ be a connected CW-complex with finite 
$1$-skeleton, with fundamental group $G=\pi_1(X)$.  
Consider the algebraic group $\T(X) =\Hom (G, \C^*)$. 
Each character $\rho\in \T(X)$ determines a rank $1$ local 
system (or, a rank $1$ complex flat bundle) on $X$, which 
we denote by $\C_{\rho}$.  The {\em characteristic varieties}\/ 
of $X$ are the jumping loci for cohomology with coefficients in 
such local systems:
\begin{equation*}
\label{eq:cv}
\VV_d(X)=\{\rho \in \T(X)
\mid \dim H^1(X,\C_{\rho})\ge d\}, \quad \text{for $d>0$}.
\end{equation*}

The characteristic varieties are Zariski closed subsets 
of $\T(X)$. These varieties depend only on the maximal 
metabelian quotient of the fundamental group, $G/G^{(2)}$, 
so we sometimes denote them as $\VV_d(G)$.  An irreducible 
component of $\VV_d(X)$ is called {\em non-translated}\/ 
if it contains the origin $1$ of the algebraic group $\T(X)$.

Consider now the cohomology algebra $H^* (X,\C)$.  
Left-mul\-tiplication by an element $x\in H=H^1(X,\C)$ 
yields a cochain complex $(H^*(X,\C), \lambda_x )$.  
The {\em resonance varieties}\/ of $X$ are the jumping 
loci for the homology of this complex:
\begin{equation*}
\label{eq:rv}
\RR_d(X)=\{x \in H \mid 
\dim H^1(H^*(X,\C),\lambda_x ) \ge  d\}, \quad \text{for $d>0$}.
\end{equation*}

The resonance varieties are homogeneous, Zariski closed subsets 
of the affine space $H^1(X, \C)= \Hom (G, \C)$. These varieties 
depend only on the co-restriction of the cup-product map, 
$\mu_G$, so we sometimes denote them by $\RR_d(G)$. 

\subsection{The exponential map}
\label{subsec:exp}

The usual exponential map, $\exp\colon \C\to \C^*$, induces 
a coefficient homomorphism, $\exp\colon H^1(G,\C) \to H^1(G,\C^*)$. 
The next result describes the behavior of this complex analytic map 
with respect to the cohomology jump loci of $G$, and some of the 
qualitative properties of these loci, under a formality assumption. 

\begin{theorem}[\cite{DPS-jump}]
\label{thm=exp}
Let $G$ be a $1$-formal group.  For each $d>0$,
\begin{enumerate}
\item \label{j1} The irreducible components of $\RR_d(G)$ 
are all linear subspaces of $H^1(G,\C)$, defined over $\Q$.
\item \label{j2}  The non-translated components of $\VV_d(G)$ 
are all subtori of the form $\exp(L)$, with $L$ running through 
the irreducible components of $\RR_d(G)$. 
\end{enumerate}
\end{theorem}

The $1$-formality hypothesis in the above theorem is crucial. 

\begin{example}
\label{ex=euclid}
Let $K$ be the finite, $2$-dimensional CW-complex defined 
in \cite[Example 4.6]{DPS-jump}.  In this case, the rationality 
property from \eqref{j1} is violated for $\RR_1(K)$. Consequently, 
the ring $H^*(K, \Q)$ cannot be realized as $H^*(X, \Q)$, for 
any {\em formal}\/ space $X$. 
\end{example}

This example stands in marked contrast with a basic result 
from simply-connected rational homotopy theory, due to 
Quillen \cite{Qu} and Sullivan \cite{Su77}:  Any finite-dimensional, 
commutative graded algebra 
$A^{*}$ defined over $\Q$, with $A^0=\Q$ and $A^1=0$, 
can be realized as the cohomology algebra of a $1$-connected, 
finite, formal CW-complex $X$.

\subsection{The BNS invariant}
\label{subsec:bns}

Let $G$ be a finitely generated group.
The following definition was introduced by Bieri, Neumann, 
and Strebel, in their seminal paper \cite{BNS}.

\begin{definition}
\label{def=bns}
Pick a finite generating set for $G$, and denote by $\CC$ the associated 
Cayley graph. The {\em BNS invariant}, $\Sigma^1(G)$, consists of those 
non-zero homomorphisms, $\varphi\colon G\to \R$, for which the full 
subgraph of $\CC$ on vertex set $\{ g\in G\mid \varphi (g)\ge 0\}$ is 
connected. The definition is independent of the choice of generators for $G$.
\end{definition}

If $M$ is a compact manifold and $p\colon M\to S^1$ is a locally trivial fibration, 
then $p^*(\omega)\in \Sigma^1(\pi_1(M))$, where $\omega\in H^1(S^1, \Z)$ 
is the generator given by the canonical orientation. In dimension $3$, the 
elements of $\Sigma^1(\pi_1(M))$ coincide with the cohomology classes in 
$H^1(M, \R)$ having closed, nowhere vanishing, de Rham representatives; 
see \cite{BNS}.

\begin{theorem}[\cite{PS-bns}]
\label{thm=bound}
If the group $G$ is $1$-formal, then 
$\Sigma^1(G) \subseteq H^1(G, \R)\setminus \RR_1(G)$.  
\end{theorem}

For right-angled Artin groups $G_{\Gamma}$, the above 
inclusion becomes equality, see \cite{PS-artin}. 

In the above theorem, the $1$-formality hypothesis is 
again crucial, as the next example illustrates.  

\begin{example}
\label{ex=nobound}
Let $M=G_{\R}/G_{\Z}$ be the $3$-dimensional Heisenberg 
nilmanifold from Example \ref{ex:triple}, with fundamental 
group $G=G_{\Z}$.  Clearly, $M$ is a torus bundle over $S^1$, 
and thus the BNS invariant $\Sigma^1(G)$ is non-empty.  On 
the other hand, $H^1(G,\R)\subseteq \RR_1(G)$, since 
$\cup_{G} =0$. Therefore, the above resonance upper 
bound for $\Sigma^1$ fails in this non-$1$-formal situation. 
\end{example}

Here is a quick application of Theorem \ref{thm=bound}. 
For a closed, orientable $3$-manifold $M$, it was shown 
in \cite{DS} that $\RR_1(M)=H^1(M,\C)$, provided $b_1(M)$ 
is even. 

\begin{corollary}[\cite{PS-bns}]
\label{cor=nofibr}
Let $M$ be a closed, orientable $3$-manifold, with even first 
Betti number. If  $M$ is $1$-formal, then $M$ does not admit  
a smooth fibration over the circle.
\end{corollary}

\section{Serre's problem and classification results}
\label{sect:serre}

We now turn to K\"{a}hler and quasi-K\"{a}hler manifolds, 
and the nature of their cohomology jumping loci. We 
illustrate the efficiency of our obstructions to $1$-formality 
and (quasi-) K\"{a}hlerianity with several classes of examples.  

\subsection{Serre's problem}
\label{subsec:serre problem}
A finitely presented group $G$ is said to be a {\em K\"{a}hler group} 
if it can be realized as $G=\pi_1(M)$, where $M$ is a compact 
K\"{a}hler manifold. If $M$ can be chosen to be a  smooth, irreducible, 
projective complex variety, then $G$ is a {\em projective group}. 
The notions of quasi-K\"{a}hler and quasi-projective group are 
defined similarly. 

Let $\mathcal{K}$, $\mathcal{P}$, $\mathcal{QK}$, 
and $\mathcal{QP}$ be the respective classes of groups.  
Clearly, $\mathcal{P}\subseteq \mathcal{K}$ and 
$\mathcal{QP}\subseteq \mathcal{QK}$, though it is 
not known whether these inclusions are strict.  Of 
course, $\mathcal{K}\subseteq \mathcal{QK}$ and 
$\mathcal{P}\subseteq \mathcal{QP}$, but both inclusions 
are strict:  for example, $\Z=\pi_1(\C^*)$ is in $\mathcal{QP}$, 
but not in $\mathcal{K}$. It is readily seen that each of these 
classes of groups is closed under finite direct products. 

\begin{problem}[J.-P.~Serre]
\label{q=serre}
Classify K\"{a}hler, projective, quasi-K\"{a}hler, and 
quasi-projec\-tive groups.
\end{problem}

This appears to be a difficult problem. As shown by Serre \cite{S}, 
all finite groups are projective. From the above discussion, it follows 
that all finitely generated abelian groups are quasi-projective. 
To the best of our knowledge, the case of nilpotent groups is open. 

\subsection{Pencils}
\label{subsec:pencils}
A good reason for grappling with Problem \ref{q=serre} is the fact 
that the fundamental group of a quasi-K\"{a}hler manifold $X$ 
determines the {\em pencils} (or, admissible maps) on $X$.

Following Arapura \cite[p.~590]{Ar}, we say that a map
$f\colon X \to C$ to a connected, smooth complex curve 
$C$ is {\em admissible}\/ if $f$ is holomorphic and surjective, 
and has a holomorphic, surjective extension with connected 
fibers to smooth compactifications, 
$\overline{f}\colon \overline{X} \to \overline{C}$, 
obtained by adding divisors with normal crossings. 
(In particular, the generic fiber of $f$ is connected, and the induced 
homomorphism, $f_{\sharp}\colon \pi_1(X)\to \pi_1(C)$, is onto.) 
Two such maps, $f\colon X \to C$ and $f'\colon X\to C'$, 
are said to be equivalent if there is an isomorphism 
$\psi\colon C \to C'$ such that $f'=\psi\circ f$.  
The pencil $f$ is called of {\em general type} if $\chi (C)<0$.
  
\begin{theorem}[Arapura \cite{Ar}]
\label{thm:arapura}
Let $X$ be a quasi-K\"{a}hler manifold, with fundamental 
group $G$.  There is a bijection between the set of positive-dimensional, 
non-translated components of $\VV_1(G)$, and the set of equivalence 
classes of pencils of general type, $f\colon X\to C$. 
This bijection associates to $f$ the component
$S_f= f^*(\T(C))$, a connected subtorus of $\T(G)$. 
\end{theorem}

The definition below, extracted from \cite{DPS-jump}, plays 
an important role in the applications.

\begin{definition}
\label{def=iso}
Let $G$ be a finitely generated group. A vector subspace 
$U\subseteq H^1(G)$ is called {\em $0$-isotropic}\/ if the 
restriction of $\cup_G$ to $U\wedge U$ is trivial. Likewise, 
$U$ is called  {\em $1$-isotropic} if this restriction is a 
non-degenerate pairing, with $1$-dimensional image.
\end{definition}

\begin{example}
\label{ex=isocurves}
Let $C$ be a complex curve with $\chi(C)<0$. Then $H^1(C)$ is 
either $1$- or $0$-isotropic, according to whether $C$ is 
compact or not.
\end{example}

Regarding the BNS invariant, here is a geometric counterpart 
to Theorem \ref{thm=bound}, based on Theorems \ref{thm:arapura} 
and \ref{thm=exp}, as well as a recent result of Delzant \cite{De1}.

\begin{theorem}[\cite{PS-bns}]
\label{thm=delbound}
Let $X$ be a compact K\"{a}hler manifold with $b_1(X)>0$, 
and let $G=\pi_1(X)$.  Then 
$\Sigma^1(G) = H^1(G, \R)\setminus \RR_1(G)$ 
if and only if there is no pencil on $X$ onto an elliptic curve, 
having multiple fibers.
\end{theorem}

\subsection{Position obstructions}
\label{subsec:pos obs}

Theorems \ref{thm=exp} and \ref{thm:arapura} lead to the 
following {\em position obstruction}, related to Serre's problem. 
Note that this obstruction depends only on $\mu_G$.

\begin{theorem}[\cite{DPS-jump}]
\label{thm=posobstr}
Let $G$ be a quasi-K\"{a}hler, $1$-formal group. Then every 
positive-dimensional irreducible component of $\RR_1(G)$ is 
$p$-isotropic with respect to $\cup_G$, and has dimension 
at least $2p+2$, for some $p\in \{0,1\}$.
\end{theorem}

This theorem  is a particular case of a more general result. 
In \cite[Theorem C]{DPS-jump}, two position obstructions 
were found, valid for an arbitrary quasi-K\"{a}hler group $G$, 
formulated in terms of the first characteristic variety $\VV_1(G)$. 
The second obstruction leads to a powerful restriction on the 
{\em multivariable Alexander polynomial}, 
$\Delta^G \in \Z [t_1^{\pm 1}, \dots, t_n^{\pm 1}]$.

\begin{theorem}[\cite{DPS07}]
\label{thm=alex}
If $G$ is a quasi-K\"{a}hler group with $n:=b_1(G)\ne 2$, 
then $\Delta^G$ has a single essential variable, that is,
$\Delta^G (t_1, \dots, t_n)= P(t_1^{e_1}\cdots t_n^{e_n})$, 
for some polynomial $P(t)\in \Z [t^{\pm 1}]$.
\end{theorem}

Theorem \ref{thm=posobstr} leads to the following complete 
classification of (quasi-) K\"{a}hler groups within the classes 
of right-angled Artin groups and Bestvina--Brady 
groups introduced in \S\S\ref{subsec:artin}--\ref{subsec:bb}. 

\begin{theorem}[\cite{DPS-jump}]
\label{thm:artinserre}
Let $\G$ be a finite simple graph, and $G_{\G}$  the 
corresponding right-angled Artin group.  Then:
\begin{enumerate}
\item $G_{\G} \in \mathcal{QK}$ $\same$ 
$G_{\G} \in \mathcal{QP}$ $\same$
$\G$ is a complete multipartite graph $K_{n_1,\dots,n_r}=
\overline{K}_{n_1} * \cdots * \overline{K}_{n_r}$.

\item $G_{\G} \in \mathcal{K}$ $\same$ 
$G_{\G} \in \mathcal{P}$ $\same$
$\G$ is a complete graph $K_{n}$, with $n$ even. 
\end{enumerate}
\end{theorem}

Note that $G_{\G}= F_{n_1}\times \cdots \times F_{n_r}$ when 
$\G =K_{n_1,\dots,n_r}$, and $G_{\G}= \Z^n$ when
$\G =K_n$. 

\begin{theorem}[\cite{DPS-sk}]
\label{thm:bbserre}
Let $\G$ be a finite simple graph, and $N_{\G}$  the 
corresponding Bestvina--Brady group. Then:
\begin{enumerate}
\item$N_{\G} \in \mathcal{QK}$ $\same$ 
$N_{\G} \in \mathcal{QP}$ $\same$ $\G$ is either a 
tree, or $\G=K_{n_1,\dots ,n_r}$, with some $n_i=1$, 
or all $n_i\ge 2$ and $r\ge 3$.

\item $N_{\G} \in \mathcal{K}$ $\same$ 
$N_{\G} \in \mathcal{P}$ $\same$
$\G=K_{n}$, with $n$ odd. 
\end{enumerate}
\end{theorem}

\subsection{Applications to $3$-manifolds}
\label{subsec:appl 3mfd}

The position obstruction from Theorem \ref{thm=posobstr} 
also turns out to be very efficient at determining which 
(quasi-) K\"{a}hler groups occur as fundamental groups 
of closed $3$-manifolds. 

\begin{theorem}[\cite{DS}]
\label{thm:3dk}
Let $G$ be the fundamental group of a closed $3$-manifold. 
Then $G$ is a K\"{a}hler group if and only if $G$ is a finite 
subgroup of $\operatorname{O}(4)$, acting freely on $S^3$.
\end{theorem}

\begin{theorem}[\cite{DPS-qk}]
\label{thm=3qk}
Let $G$ be the fundamental group of a closed, orientable 
$3$-manifold.  Then (up to Malcev completion) $G$ is a 
quasi-K\"{a}hler, $1$-formal group if and only if
$G=F_n$, or $G=\Z\times \pi_1(\Sigma_g)$.
\end{theorem}

In the case of boundary manifolds of line arrangements 
in $\CP^2$, more can be said.  Let $\A$ be such an 
arrangement, and let $M$ be the closed, orientable 
$3$-manifold obtained by taking the boundary of 
a regular neighborhood of $\A$ in $\CP^2$. 
Theorem \ref{thm=exp} was used in \cite[Theorem 9.7]{CS08} 
to classify those boundary manifolds which are formal,   
while Theorem \ref{thm=alex} was used 
in \cite[Proposition 4.7]{DPS07} 
to classify those boundary manifolds whose 
fundamental groups are quasi-projective.  
We summarize these results, as follows. 

\begin{theorem}[\cite{CS08, DPS07}] 
\label{thm:bdry}
Let $\A=\{\ell_0,\dots,\ell_n\}$ be an arrangement of lines in 
$\CP^2$, and let $M$ be the corresponding boundary manifold. 
The following are equivalent:
\begin{enumerate}
\item  \label{f1} The manifold $M$ is formal. 
\item  \label{f2} The group $G=\pi_1(M)$ is $1$-formal.
\item  \label{f3} The group $G$ is quasi-projective.
\item  \label{f4} $\A$ is either a pencil or a near-pencil. 
\end{enumerate}
\end{theorem}

The corresponding $3$-manifolds are easy to describe: if 
$\A$ is a pencil, then $M=\sharp^n S^1\times S^2$, and if  
$\A$ is a near-pencil, then $M=S^1\times \Sigma_{n-1}$.

\section{Formality and the monodromy action}
\label{sect:mono}

We go on by examining the interplay between algebraic 
monodromy and $1$-formality, with emphasis on extensions 
of $\Z$, and fibrations over $S^1$. In this section, all manifolds 
are compact and connected.

\subsection{Artin kernels}
\label{ss61}

Let $\G$ be a finite simple graph. Every epimorphism 
$\chi\colon G_{\G}\surj \Z$ from the right-angled Artin group 
$G_{\G}$ to the integers gives rise to an {\em Artin kernel}, 
$N_{\chi}= \ker (\chi)$,  generalizing the Bestvina--Brady
group $N_{\G}=\ker(\nu)$ from \S\ref{subsec:bb}.  

In \cite{PS-toric}, we found a combinatorial procedure which characterizes 
the triviality of the monodromy action of $\Z$ on the homology groups 
of $N_{\chi}$, up to a fixed degree $q$. As a result, we were able to 
establish the $1$-formality of a large class of Artin kernels. 

\begin{theorem}[\cite{PS-toric}]
\label{thm:artink}
If the monodromy action on $H_*(N_{\chi})$ is trivial, up to 
degree $2$, then the Artin kernel $N_{\chi}$ is a $1$-formal group.
\end{theorem}

For the Bestvina--Brady groups $N_{\G}$, this triviality test boils 
down to verifying that $\widetilde{H}_i (\Delta_{\G})=0$, for $i\le 1$.  
In the case when the flag complex $\Delta_{\G}$ is simply-connected, 
we recover Theorem \ref{thm:formalbb}.  

\begin{remark}
\label{rem:mono heis}
The above theorem may be regarded as a complement 
to \S\ref{subsec:operations}, indicating a new formality-preserving 
construction: the passage from $1$-formal groups to normal 
subgroups with infinite cyclic quotient, under certain triviality 
assumptions on the monodromy action.  That the assumption 
on the monodromy is key to preserving formality is illustrated 
by the $1$-formal group $H=G\rtimes_{\varphi} \Z$ from 
Example \ref{ex=subgr}. In that case, we have 
an epimorphism $\chi\colon H\surj \Z$, with kernel the 
(non-$1$-formal) Heisenberg group $G$, 
and non-trivial monodromy action on $H_1(G)$. 
\end{remark}

Theorem \ref{thm:artink} enabled us to construct what seem to be 
the first instances of $1$-formal groups which are {\em not}\/ 
finitely presented. The example below is taken 
from \cite[Example 10.3]{PS-toric}. 

\begin{example}
\label{ex=fnotfp}
Let $L=\Delta_\Gamma$ be a flag triangulation of the real 
projective plane, $\RP^2$.  Clearly, $L$ is connected. 
On the other hand, $H_1(L,\Z)=\Z_2$, and so, by \cite{BB}, 
$N_{\Gamma}$ is not finitely presented.  But $H_1(L,\Q)=0$, 
and so, by Theorem \ref{thm:artink}, $N_{\Gamma}$ is $1$-formal. 
\end{example}

\subsection{Jordan blocks}
\label{ss62}

We examine now the converse question: does $1$-formality impose 
restrictions on the algebraic monodromy? In \cite[Proposition 9.4]{PS-spectral}, 
we obtained a general result, relating the monodromy action in 
extensions of $\Z$ to the resonance varieties, without any 
formality assumptions.  We used this to deduce in \cite{PS-mono} 
the following implication of $1$-formality on the algebraic monodromy.

Let $U$ be a closed manifold, and let $h\colon U\to U$ be a
diffeomorphism.   Denote by $U_h$ the mapping torus of $h$. 

\begin{theorem}[\cite{PS-mono}] 
\label{thm:fgm}
Let $p\colon M\to U_{h}$ be a locally trivial smooth fibration. 
Assume $M$ is $1$-formal. Then, if the monodromy operator 
$h_* \colon H_1(U)\to H_1(U)$ has eigenvalue $1$, all corresponding 
Jordan blocks must have size $1$.
\end{theorem}

This theorem substantially extends a result of Fern\'{a}ndez, Gray, and 
Morgan \cite{FGM}, valid only for circle bundles over $U_h$. Those 
authors proved their result by a different method, which relies on 
Massey products.

\begin{example}
\label{ex:heisagain}
The $1$-formality hypothesis is crucial here, even in the 
particular case when $M=U_h$ and $p=\id$.  For 
instance, the Heisenberg manifold from Example \ref{ex:triple} is not 
$1$-formal, yet it fibers over $S^1$, with fiber the $2$-torus and monodromy 
$\left( \begin{smallmatrix} 1 & 1\\ 0 & 1 \end{smallmatrix} \right)$.  
\end{example}

\begin{example}
\label{ex:sing}
Let $f\colon (\C^2, 0)\to (\C, 0)$ be the germ of a reduced polynomial 
function, with associated singularity link $K$, and Milnor fibration 
$p\colon S^3\setminus K\to S^1$. A basic fact in singularity theory 
is that the algebraic monodromy of $p$ has no Jordan block
of size greater than $1$, for the eigenvalue $1$. As shown by 
Durfee and Hain \cite{DH}, the link complement is a formal space. 
Hence, the result on the algebraic monodromy follows from 
Theorem \ref{thm:fgm}, applied to the link exterior.
\end{example}

\section{Epilogue: beyond formality}
\label{sect:epilogue}

We conclude with a class of examples inspired by the 
work of Geiges on $2$-torus bundles over the $2$-torus \cite{Ge}.   
These examples illustrate some of the subtle interplay 
between monodromy and geometric structures on 
manifolds.  

\subsection{Mapping tori and K\"{a}hler metrics}
\label{subsec:map tori}

Let $U=\Sigma_g$ be a compact Riemann surface of 
genus $g\ge 1$, and let $h\colon U \to U$ be an 
orientation-preserving diffeomorphism (which may be 
viewed as an element of the mapping class group 
$\mathcal{M}_g$). 

As before, denote by $U_h$ the mapping torus of $h$. 
By construction, this is a closed, orientable $3$-manifold 
fibering over the circle, with fiber $U$ and monodromy $h$. 

\begin{question}
\label{question:geiges}
Suppose $1$ is not an eigenvalue of the algebraic 
monodromy operator, $h_* \colon H_1(U)\to H_1(U)$.   
Does there exist a closed, connected manifold $N$ 
such that $N\times U_h$ carries a K\"{a}hler metric?
\end{question}

In general, a question of this type has a negative answer.  Using a deep 
classification result of Wall \cite{W} concerning complex structures 
on $4$-manifolds, Geiges showed in \cite[p.~555]{Ge} that the 
manifold $S^1 \times (T^2)_h$, where $h=\left( \begin{smallmatrix} 
-1 & 0\\ 0& -1 \end{smallmatrix} \right)$, supports no K\"{a}hler metric. 

On the other hand, a negative answer to 
Question \ref{question:geiges} cannot be obtained solely 
by rational homotopy methods.  The next result makes it clear why. 

\begin{prop}
\label{prop:rht geiges}
With the setup from Question \ref{question:geiges}, 
suppose $N=S^1\times M$, where $M$ is a compact K\"{a}hler 
manifold.  Then $N \times U_h$ has the same $\Q$-homotopy 
type as a K\"{a}hler manifold, namely, $M\times T^2\times S^2$. 
\end{prop}

\begin{proof}
Recall we are assuming $1$ is not an eigenvalue of 
$h_* \colon H_1(U)\to H_1(U)$, that is, $\id - h_*$ is an 
isomorphism.  From the Wang sequence of the fibration 
$U\to U_h \to S^1$, we infer that $b_1(U_h)=1$.  
Now, since $U_h$ is a closed, orientable $3$-manifold, 
Theorem \ref{thm:fm} implies that $U_h$ is formal, 
and has the same $\Q$-homotopy type as $S^1\times S^2$.  
Our claim follows at once. 
\end{proof}

\subsection{Eigenvalues of the monodromy operator}
\label{subsec:evalues}

For a space $X$, recall $\T(X)$ denotes the character group 
$\Hom(\pi_1(X),\C^*)$.  Let $\T^0(X)$ be the connected component 
of $1\in \T(X)$.  

Returning to the setting of Question \ref{question:geiges}, 
let $p\colon U_h \to S^1$ be the canonical fibration of the 
mapping torus. Since $b_1(U_h)=1$, the induced homomorphism, 
$p^*\colon H^1(S^1,\C^*) \to H^1(U_h,\C^*)$, may be
used to identify $\T^0(U_h)$ with $\T(S^1)=\C^*$.

\begin{lemma}
\label{lem:eigen}
Under the above identification, the variety 
$\VV_1(U_h) \cap \T^0(U_h)\subset \C^*$ 
consists of $1$, together with the eigenvalues of $h_*$.
\end{lemma}

\begin{proof}
Fix a character $\rho\in \T(S^1)=\C^*$.  
Consider the Leray-Serre spectral sequence of the fibration 
$U \inj U_h \xrightarrow{\,p\,} S^1$, with coefficients in 
the local system determined by $\rho$: 
\[ 
E^2_{s,t}=H_s(S^1,H_t(U,\C)) \Rightarrow 
H_{s+t}(U_h, \C_{p^*(\rho)}),
\]
where $\pi_1(S^1)$ acts on $H_t(U, \C)$ by $\rho^{-1}\cdot h_*$. 
In total degree $1$, we obtain an isomorphism 
\[
H_1 (U_h, \C_{p^*(\rho)})\cong 
H_1 (S^1, \C_{\rho}) \oplus \coker (h_* -\rho \cdot \id) .
\]
Now, $H_1 (S^1, \C_{\rho})=0$ or $\C$, according to 
whether $\rho\ne 1$ or $\rho =1$.  The conclusion follows. 
\end{proof}

\subsection{Quasi-K\"{a}hler manifolds and monodromy}
\label{subsec:qkmono}

For a quasi-K\"{a}hler manifold $X$, a powerful result, due 
to Arapura \cite{Ar}, guarantees that the isolated points of 
$\VV_1(X)$ must be unitary characters.  Within the realm 
of groups we consider here, this leads to a simple-to-verify 
obstruction for membership in $\mathcal{QK}$.  

\begin{prop}
\label{prop=unit}
Let $N$ be a compact, connected manifold, and 
suppose $\pi_1(N\times U_h)$ is a quasi-K\"{a}hler group. 
Then all eigenvalues of the monodromy operator, 
$h_*\colon H_1(U)\to H_1(U)$, have norm $1$. 
\end{prop}

\begin{proof}
By the K\"unneth formula, 
\[
\VV_1(N\times U_h)= 
1\times \VV_1(U_h) \cup \VV_1(N)\times 1.
\]
By Lemma \ref{lem:eigen}, each eigenvalue $\rho$ of $h_*$ 
gives rise to the isolated point $1\times \rho$ of $\VV_1(N\times U_h)$.  
In view of the aforementioned result of Arapura, this finishes 
the proof.
\end{proof}

It is an easy matter to construct elements $h\in \mathcal{M}_g$ 
such that $h_*$ has no eigenvalue of norm $1$.  

\begin{example}
\label{ex:nonunit}
Pick a matrix $A \in \SL(2, \Z)$ with $\abs{\tr(A)} \ge 3$.  
Then $A$ has two distinct, non-unitary eigenvalues, say, 
$\lambda_1$ and $\lambda_2$.  Let $B$ be the block-sum 
of $g$ copies of $A$. Clearly, $B$ belongs to $\Sp(2g, \Z)$, 
and has the same eigenvalues as $A$. By a classical result, 
there exists a diffeomorphism $h\colon U\to U$ (necessarily, 
orientation-preserving), such that $h_*\colon H_1(U) \to H_1(U)$ 
has matrix $B$. 
\end{example}

Proposition \ref{prop=unit} may be applied to give 
a negative answer to Question \ref{question:geiges}, even 
when rational homotopy methods are inconclusive.

\begin{corollary}
\label{cor:nonkahler}
Let $U=\Sigma_g$ and let $h\colon U\to U$ be a diffeomorphism 
as constructed in Example \ref{ex:nonunit}.  Then, for any closed 
manifold $N$, the group $\pi_1(N\times U_h)$ is not a quasi-K\"{a}hler group. 
In particular, the manifold $N\times U_h$ does not carry any 
K\"{a}hler metric. 
\end{corollary}

Putting things together, we can now prove the result stated 
in the Introduction. 

\begin{proof}[Proof of Theorem 1.1]
With notation as in the previous corollary, let $W=S^1\times U_h$. 
Clearly, $W$ is a closed, orientable, formal $4$-manifold.  
Let $M$ be an arbitrary compact K\"{a}hler manifold, and 
write $N=S^1\times M$.  By Proposition \ref{prop:rht geiges}, 
$M\times W=N\times U_h$ has the $\Q$-homotopy 
type of a K\"{a}hler manifold, namely, $M\times T^2\times S^2$. 
On the other hand, by Corollary \ref{cor:nonkahler},  
$M\times W$ admits no K\"{a}hler metric. 

It remains to produce infinitely many manifolds $W$ as above.  
For that, choose a surface $U=\Sigma_g$ and matrices 
$A_n=\left( \begin{smallmatrix} n+2 & -1\\ 1 & 0 \end{smallmatrix} \right)$, 
with $n> 1$, as input for the construction. Denote by $B_{g,n}$ 
the block-sum of $g$ copies of $A_n$, and let $W_{g,n}$ be 
the resulting $4$-manifold. 
A straightforward computation shows that 
$H_1(W_{g,n}, \Z)=\Z^2\oplus \bigoplus^g \Z/n\Z$.  
Hence, the manifolds $\{ W_{g,n}\}_{g\ge 1,\, n>1}$ 
are pairwise distinct.
\end{proof} 

\newcommand{\arxiv}[1]
{\texttt{\href{http://arxiv.org/abs/#1}{arXiv:#1}}}
\renewcommand{\MR}[1]
{\href{http://www.ams.org/mathscinet-getitem?mr=#1}{MR#1}}

\end{document}